\documentclass[reqno,12pt,a4paper]{amsart}
\pdfoutput=1
\usepackage[ascii]{inputenc}
\usepackage[bookmarksnumbered,pdfusetitle]{hyperref}
\usepackage{bbm,microtype,amsmath,amsthm,amssymb,mathtools,graphicx,appendix,braket}
\usepackage[usenames,dvipsnames,svgnames,table]{xcolor}
\usepackage{amsrefs}
\usepackage[capitalize]{cleveref}
\newtheorem{theorem}{Theorem}[section]
\newtheorem{proposition}[theorem]{Proposition}

\newtheorem{definition}[theorem]{Definition}

\newtheorem{remark}[theorem]{Remark}
\newtheorem{example}[theorem]{Example}
\numberwithin{equation}{section}
\DeclareMathOperator{\coker}{coker}
\DeclareMathOperator{\Cone}{Cone}
\DeclareMathOperator{\Ext}{Ext}

\DeclareMathOperator{\eul}{eul}
\DeclareMathOperator{\GL}{GL}

\DeclareMathOperator{\Gr}{Gr}
\DeclareMathOperator{\Horn}{Horn}
\DeclareMathOperator{\Hom}{Hom}
\DeclareMathOperator{\Sym}{Sym}
\DeclareMathOperator{\ext}{ext}
\DeclareMathOperator{\linspan}{span}
\newcommand{\C}{{\mathbb C}}
\newcommand{\N}{{\mathbb N}}

\newcommand{\ZZ}{{\mathbb Z}}

\newcommand{\CF}{{\mathcal F}}
\newcommand{\CG}{{\mathcal G}}
\newcommand{\CH}{{\mathcal H}}

\newcommand{\CJ}{{\mathcal J}}
\newcommand{\CK}{{\mathcal K}}
\newcommand{\CL}{{\mathcal L}}

\newcommand{\CN}{{\mathcal N}}

\newcommand{\CS}{{\mathcal S}}

\newcommand{\CU}{{\mathcal U}}
\newcommand{\CV}{{\mathcal V}}

\newcommand{\CW}{{\mathcal W}}

\newcommand{\g}{{\mathfrak{g}}}
\newcommand{\gl}{{\mathfrak{gl}}}

\newcommand{\n}{{\mathfrak{n}}}

\renewcommand{\t}{{\mathfrak{t}}}

\allowdisplaybreaks[4]
\begin{document}
\title{Horn inequalities and quivers}
\author{V.~Baldoni}
\address{Velleda Baldoni: Dipartimento di Matematica, Universit\`a degli studi di Roma ``Tor Vergata'', Via della ricerca scientifica 1, I-00133, Italy}
\email{baldoni@mat.uniroma2.it}
\author{M.~Vergne}
\address{Mich\`ele Vergne: Universit\'e Paris 7 Diderot, Institut Math\'ematique de Jussieu, Sophie Germain, case 75205, Paris Cedex 13}
\email{vergne@math.jussieu.fr}
\author{M.~Walter}
\address{Michael Walter: QuSoft, Korteweg-de Vries Institute for Mathematics, Institute for Theoretical Physics, Institute for Logic, Language and Computation, University of Amsterdam, 1098 XG Amsterdam, Netherlands}
\email{m.walter@uva.nl}
\date{}
\hypersetup{pdfauthor={V. Baldoni, M. Vergne, M. Walter}}


\begin{abstract}
Let $G$ be a complex reductive group acting on a finite-dimensional complex vector space~$\CH$.
Let $B$ be a Borel subgroup of~$G$ and let $T$ be the associated torus.
The Mumford cone is the polyhedral cone generated by the $T$-weights of the polynomial functions on~$\CH$ which are semi-invariant under the Borel subgroup.
In this article, we determine the inequalities of the Mumford cone in the case of the linear representation of the group $G=\prod_x \GL(n_x)$ associated to a quiver and a dimension vector $n=(n_x)$.
We give these inequalities in terms of filtered dimension vectors, and we directly adapt Schofield's argument to inductively determine the dimension vectors of general subrepresentations in the filtered context.
In particular, this gives one further proof of the Horn inequalities for tensor products.
\end{abstract}

\maketitle

\section{A generalized Horn condition}\label{sec:general}
Let $G$ be a connected complex reductive algebraic group, $K$ a maximal compact subgroup of~$G$ and $T$ be a maximal torus of~$K$.
Let $\t$ denote the Lie algebra of~$T$, and $\Lambda\subset i\t^*$ be the lattice of weights of~$T$.
Let us choose a Borel subgroup~$B$ of~$G$ which contains~$T$.
We denote the set of associated positive roots by $R^+_G$.
This determines a closed Weyl chamber~$i\t^*_{\geq 0}$ and the set~$P_G$ of dominant weights.
We denote by~$V_\lambda$ the irreducible representation of $G$ associated to $\lambda\in P_G$.

Now let $\Pi \colon G \rightarrow \GL(\CH)$ be a linear representation of $G$ on a finite-dimensional complex vector space $\CH$.
We denote by~$\pi$ the corresponding infinitesimal representation of the Lie algebra of~$G$.
We assume that the action of $G$ in the space $R(\CH)$ of polynomial functions on~$\CH$ decomposes with finite multiplicities.
This is the case, e.g., when the action includes the homotheties.
Consider the \emph{Mumford cone} $\Cone(G,\CH)\subset i\t^*$, defined as the cone generated by the dominant weights $\lambda$ such that $V_{\lambda} \subseteq R(\CH)$.
The equation $V_{\lambda} \subseteq R(\CH)$ means that there exists a nonzero polynomial function on $\CH$, which is semi-invariant under the action of $B$ and of weight $\lambda$ under $T$.
The cone $\Cone(G,\CH)$ is a polyhedral cone determined by a finite number of inequalities.
By Mumford's theorem~\cite{NessMumford84}, it can also be described using the moment map.

\begin{example}\label{exa:Horn}
Consider $\CH=\bigoplus_{i=1}^s \Hom(V_{x_i}, V_y)$, where $V_{x_i}$, $V_y$ are $n$-di\-men\-si\-o\-nal vector spaces.
It has a natural action of $G=\prod_{i=1}^s \GL(n) \times \GL(n)$.
The cone $\Cone(G,\CH)$ is generated by the tuples of dominant weights $(\lambda_1,\dots,\lambda_s,\mu)$ associated to Young diagrams such that $V_\mu \subseteq \bigotimes_{i=1}^s V_{\lambda_i}.$
It is described by the well-known Horn inequalities~\cites{MR0140521,MR1671451}.
The $G$-representation $\CH$ is naturally associated to the quiver $x_i \to y$, $i=1,\dots,s$, and the dimension vector $n_{x_i}=n_y=n$ (see \cref{eq:quiver repr,sec:q-intersection} for precise definitions).
\end{example}

Let $\Phi_\CH \subset i\t^*$ be the list of weights for the action of $T$ in $\CH$.
When $G=T$ is a torus, $\Cone(T,\CH)$ is the cone generated by the elements~$-\phi$ for $\phi\in \Phi_{\CH}$.

\begin{definition}\label{def:admissible}
  An element $H \in i\t$ is called \emph{admissible} if the linear hyperplane $\{H=0\}$ is spanned by a subset of $\Phi_\CH$.
\end{definition}

Let $H \in i\t$. Decompose
\begin{equation}
\label{eq:weight space decomposition}
  \CH = \CH(H < 0) \oplus \CH(H = 0) \oplus \CH(H > 0),
\end{equation}
where $\CH(H < 0) = \bigoplus_{\phi\in \Phi_\CH : (\phi,H) < 0} \CH_\phi$ is the sum of the eigenspaces of the Hermitian operator $\pi(H)$ with eigenvalue less than $0$, etc.

If all weights $\phi\in \Phi_\CH$ are contained in the halfspace $\{H\geq0\}$ associated with the hyperplane $\{H=0\}$, it is immediate to see that the cone $\Cone(G,\CH)$ lies in the opposite half-space $\{H\leq0\}$.
This characterizes the Kirwan cone when $G=T$ is a torus.

However, in general, there are facets of $\Cone(G,\CH)$ determined by some admissible element~$H$ with weights on \emph{both} sides of the hyperplane~$\{H=0\}$.
In this case, they are ``balanced'' by the roots $R^+_G$ of~$G$ in a sense that we describe now (see \cite{VW}).
Let $\mathfrak n$ be the Lie algebra of the maximal unipotent subgroup~$N$ of $B$.
Similarly to above, let~$\n(H<0)$ denote the sum of all root spaces for positive roots $\alpha$ such that $(\alpha,H)<0$.
Lastly, consider for $v\in \CH(H=0)$ the map
\begin{align*}
  \delta_v \colon \n(H < 0) \rightarrow \CH(H < 0), \quad X \mapsto \pi(X)v.
\end{align*}

\begin{definition}\label{def:euler mumford}
We define
\begin{align*}
  \eul(\CH,H)= \dim(\n(H < 0))-\dim \CH(H < 0)
\end{align*}
the difference of the dimensions of the source space and the target space of the map~$\delta_v$.
Furthermore, define $\kappa_H \in i\t^*$ to be the difference of the traces of the action of $i\t$ in the source space and the target space.
That is, for $X\in i\t$,
\begin{align*}
  \kappa_H(X) \coloneqq \!\!\!\!\sum_{\alpha \in R^+_G : (\alpha,H) < 0} \!\!\!\!\alpha(X)-\sum_{\phi \in \Phi_\CH : (\phi,H) < 0} \!\!\!\!\phi(X).
\end{align*}
\end{definition}

If $\eul(\CH,H)=0$ then we can consider the polynomial function $P_H$ on $\CH(H=0)$ defined by
\begin{align*}
  P_H(v) \coloneqq \det(\delta_v)
\end{align*}
where we evaluate the determinant with respect to arbitrary fixed bases of the source and target spaces.
This polynomial function is nonzero (i.e., not everywhere vanishing) if and only if the map $\delta_v$ is an isomorphism for generic $v$ in $\CH(H=0)$.

\begin{definition}\label{def:covering ressayre}
We shall say that $H\in i\t$ is a \emph{covering element} if
\begin{enumerate}
  \item $H$ is admissible and
  \item the map $\delta_v$ is surjective for generic~$v$.
\end{enumerate}
We shall say that $H\in i\t$ is a \emph{Ressayre element} if
\begin{enumerate}
  \item $H$ is admissible,
  \item $\eul(\CH,H)=0$, and
  \item the function $P_H$ is a nonzero polynomial function.
\end{enumerate}
\end{definition}

Clearly, every Ressayre element is a covering element.
If $H$ is a covering element, it is easy to prove that the cone $\Cone(G,\CH)$ is contained in the half-space $\{H\leq0\}$.
See~\cite{VW}, but beware of a change of notation since there we were considering the cone generated by the dominant weights $\lambda$ such that $V_\lambda\subseteq R(\CH)^* = \Sym(\CH)$.
Thus, covering elements and in particular Ressayre elements yield inequalities on $\Cone(G,\CH)$ (which may or may not be redundant).
It follows from Ressayre's results~\cites{ressayre2010geometric,ressayre2011geometric} that these inequalities characterize the Mumford cone (see \cite{VW} for a simple proof in the case the cone is solid).

We now come to a central definition of this paper.

\begin{definition}\label{def:horn}
We say that $H\in i\t$ is a \emph{Horn element} if
\begin{enumerate}
\item $H$ is admissible,
\item $\eul(\CH,H)=0$, and
\item $\kappa_H \in \Cone(G(H=0), \CH(H=0))$.
\end{enumerate}
\end{definition}

Thus, we replace the condition that $P_H\neq0$ for a Ressayre element by an inductive condition in terms of a moment polytope for a smaller group and representation.
The group~$G(H=0)$ is the centralizer of the element~$H$.
The Borel subgroup $B$ of $G$ determines a Borel subgroup~$B(H=0)$ with the same torus~$T$. 

Every Ressayre element is also a Horn element.
Indeed, $P_H$ is a polynomial on $\CH(H=0)$ that is semi-invariant by~$B(H=0)$, with weight~$\kappa_H$.
Thus, if $P_H$ is nonzero then $\kappa_H$ is an element of the Mumford cone associated with the action of $G(H=0)$ on $\CH(H=0)$.

In \cref{exa:Horn}, it is easy to see that the inequalities $\{H\leq0\}$ determined by the Horn elements~$H$ are exactly the inequalities conjectured by Horn~\cite{MR0140521} and proved by Knutson-Tao~\cite{MR1671451}.
Thus they are a complete set of inequalities that describe $\Cone(G,\CH)$.
In view of this emblematic example, we give the following definition.

\begin{definition}\label{def:horn prop}
A representation $\Pi$ of $G$ in $\CH$ has the \emph{Horn property} if the Mumford cone is defined by
\begin{align*}
  \Cone(G,\CH) = \left\{ \lambda \in i \t^*_{\geq 0} : (H, \lambda) \leq 0 \text{ for all Horn elements $H$} \right\}.
\end{align*}
\end{definition}

When the representations $\CH(H=0)$ again satisfy the Horn property and so forth (as in the case of \cref{exa:Horn}), one obtains an inductive system of inequalities describing the moment polytope.

\medskip

We now describe our main result.
Let $Q=(Q_0,Q_1)$ be a quiver without cycles.
Consider a family $\CV$ of complex vector spaces $V_x$ for~$x\in Q_0$ and the space
\begin{align}\label{eq:quiver repr}
  \CH_Q(\CV) = \bigoplus_{a : x\to y\in Q_1} \Hom(V_x, V_y).
\end{align}
We have a natural action of $\GL_Q(\CV) = \prod_{x\in Q_0} \GL(V_x)$.
Our main result is that the action of $\GL_Q(\CV)$ on $\CH_Q(\CV)$ has the Horn property (\cref{thm:main}, in fact, a weaker inductive ``half condition'' characterizes the moment cone).

Previously, Derksen-Weyman~\cite{MR1758750} and Schofield-Van den Bergh~\cite{MR1908144} had described the inequalities of the cone~$\Sigma_Q(\CV)$ generated by weights of the polynomials on $\CH_Q(\CV)$ that are semi-invariant \emph{under the full group} $\GL_Q(\CV)$, in terms of dimension vectors of general subrepresentations.
Furthermore, Schofield~\cite{MR1162487} gave a recursive numerical criterium to determine these dimension vectors.
Using an augmented quiver~$\tilde Q$, Derksen-Weyman deduced the Horn inequalities in \cref{exa:Horn} from their description of the corresponding semi-invariant cone~$\Sigma_{\tilde Q}(\tilde \CV)$.

We employ a more direct approach, partly inspired by Derksen-Weyman and partly by our previous study of inequalities of general moment cones~\cite{VW}.
We prove directly analogs of Schofield's results~\cite{MR1162487} on dimension vectors, for ``filtered'' dimension vectors and $\Ext$ groups in the context of vector spaces $V_x$ with filtrations (see \cref{sec:q-intersection,sec:ext}).
This allows us to prove the Horn property for arbitrary quivers and, in particular, yields one more proof of the Horn inequalities.
The description of the inequalities for the cone $\Sigma_Q(\CV)$ given by Derksen-Weyman and Schofield-Van den Bergh can also be deduced as a special case.

Our calculations of the $\Ext$ groups for filtered dimension vectors are very similar to Belkale's proof of the necessary and sufficient conditions for Schubert cells to intersect~\cite{MR2177198} (see also~\cite{BVW}).
They will be presented in full detail in a forthcoming article.

\section{\texorpdfstring{Filtered dimension vectors and $Q$-intersection}{Filtered dimension vectors and Q-intersection}}\label{sec:q-intersection}
Let $Q=(Q_0,Q_1)$ be a quiver without cycles, where $Q_0$ is the set of vertices and $Q_1$ the set of arrows.
A \emph{dimension vector} is a family of nonnegative integers~$n=(n_x)_{x\in Q_0}$.
To any family of vectors spaces $\CV=(V_x)_{x\in Q_0}$ we can associate a dimension vector by $n_x=\dim V_x$.

A \emph{representation} of $Q$ consists of a collection of vector spaces~$\CV$ as above, together with a collection of maps $v=(v_a)_{a\in Q_1} \in \CH_Q(\CV)$, indexed by the arrows $a$ of $Q$.
Let $\CS=(S_x)_{x\in Q_0}$ be a collection of subspaces $S_x \subseteq V_x$.
Then $\CS$ is called a \emph{subrepresentation} of $(\CV,v)$ if $v_a S_x\subseteq S_y$ for every arrow $a: x\to y$ in $Q_1$.

Schofield~\cite{MR1162487} determined by induction the dimension vectors~$\alpha$ such that for every $v \in \CH_Q(\CV)$ there exists a subrepresentation~$\CS$ of $(\CV,v)$ such that $\dim S_x=\alpha_x$.

\begin{definition}\label{def:filtration}
  A (complete) \emph{filtration} $F$ on a vector space $V$ is a chain of subspaces
  \[ \{0\} = F(0) \subseteq F(1) \subseteq \cdots \subseteq F(i) \subseteq F({i+1}) \subseteq \dots \subseteq F(\ell) = V, \]
  such that $\dim F(i+1)\leq\dim F(i)+1$ for all $i=0,\dots,\ell-1$, i.e., the dimensions increase by at most one in each step.
\end{definition}

The distinct subspaces in a filtration determines a flag.
However, note that the subspaces $F(i)$ need \emph{not} be strictly increasing.
Thus, if $S$ is a subspace of $V$, then $S$ inherits the filtration $F_S(i) \coloneqq F(i)\cap S$, and the quotients space $V/S$ inherits a filtration $F_{V/S}(i)\coloneqq(F(i)+S)/S$.

\begin{definition}
A \emph{filtered dimension vector} is a pair $(\CV,\CF)$, where $\CV=(V_x)_{x\in Q_0}$ is a family of vector spaces and $\CF=(F_x)_{x\in Q_0}$ is a family of filtrations, with each $F_x$ a filtration on $V_x$.
\end{definition}

Any filtered dimension vector $(\CV,\CF)$ determines a Borel subgroup which we shall denote by $B_Q(\CV,\CF) \subseteq \GL_Q(\CV)$.

Now let $\CS\subseteq\CV$ be a family of subspaces $S_x \subseteq V_x$, with dimension vector $\alpha_x=\dim S_x$.
Thus, $\CS$ determines a point in
\begin{align*}
  \Gr(\alpha,\CV) \coloneqq \prod_x \Gr(\alpha_x, V_x),
\end{align*}
where $\Gr(\alpha_x,V_x)$ denotes the Grassmannian of subspaces of dimension~$\alpha_x$ of~$V_x$.
The orbit of $\CS$ by $B_Q(\CV,\CF)$ is a product of Schubert cells in $\Gr(\alpha_x,V_x)$.
We denote the closure of this orbit by $\Omega(\CS,\CF)$ and say that it is a \emph{Schubert variety}.

For the following definition, fix a filtered dimension vector $(\CV,\CF)$.

\begin{definition}\label{def:q-intersecting}
Let $\alpha$ be a dimension vector and $\Omega \subseteq \Gr(\alpha,\CV)$ a Schubert variety.
We say that $\Omega$ is \emph{$Q$-intersecting} if, for every $v\in\CH_Q(\CV)$, there exists a subrepresentation $\CN$ of $(\CV,v)$ such that $\CN\in\Omega$.
We write~$\CS\subseteq_Q \CV$ if the Schubert variety~$\Omega(\CS,\CF)$ is $Q$-intersecting.
\end{definition}

In the case of the quiver $x_i\to y$, $i=1,\ldots,s$ and $(\CV,\CF)$ with $\dim V_{x_i}=\dim V_y=n$, which corresponds to \cref{exa:Horn}, the problem of determining the $Q$-intersecting Schubert varieties with dimension vector $\alpha_{x_i}=\alpha_y=r$ is the same problem as determining the intersecting Schubert classes in $\Gr(r,n)$ (see, e.g., \cite{BVW}).

\section{\texorpdfstring{$\Ext$ groups for filtered dimension vectors}{Ext groups for filtered dimension vectors}}\label{sec:ext}
For two families $\CV=(V_x)_{x\in Q_0}$, $\CW=(W_x)_{x\in Q_0}$ of vector spaces, define
\begin{align*}
  \CH_Q(\CV,\CW) &\coloneqq \!\!\!\!\!\bigoplus_{a:x\to y \in Q_1}\!\!\!\!\! \Hom(V_x,W_y), \\
  \g_Q(\CV,\CW) &\coloneqq \bigoplus_{x\in Q_0} \Hom(V_x,W_x).
\end{align*}
If $\CV=\CW$, the space $\CH_Q(\CV,\CV)$ is simply $\CH_Q(\CV)$ as introduced in~\eqref{eq:quiver repr}, and $\gl_Q(\CV)\coloneqq\g_Q(\CV,\CV)$ is the Lie algebra of $\GL_Q(\CV)$.

Thus, if $\alpha$ is the dimension vector of~$\CV$ and $\beta$ the dimension vector of~$\CW$, the \emph{Euler form}
\begin{align*}
  \braket{\alpha,\beta} = \sum_{x\in Q_0} \alpha_x\beta_x - \!\!\!\!\sum_{a:x\to y\in Q_1} \alpha_x\beta_y
\end{align*}
is precisely equal to $\dim \g_Q(\CV,\CW)-\dim\CH_Q(\CV,\CW)$.

Let $\Phi=(\Phi_x)_{x\in Q_0}$ be an element of $\g_Q(\CV,\CW)$, $v=(v_a)_{a\in Q_1}$ in~$\CH_Q(\CV)$, and $w=(w_a)_{a\in Q_1}$ in~$\CH_Q(\CW)$.
We denote by $\Phi v-w\Phi$ the element of $\CH_Q(\CV,\CW)$ that assigns
\begin{align*}
  \Phi_y v_a - w_a \Phi_x \colon V_x \to W_y.
\end{align*}
to each arrow $a:x\to y \in Q_1$.

\begin{definition}\label{def:gQ}
Let $(\CV,\CF)$, $(\CW,\CG)$ be filtered dimension vectors.
We say that $\Phi=(\Phi_x)\in \g_Q(\CV,\CW)$ is \emph{compatible} with the filtrations $\CF$, $\CG$ if, for all $x\in Q_0$ and $i$,
\begin{align*}
  \Phi_x(F_x(i))\subseteq G_x(i).
\end{align*}
We denote by $\g_{Q,\CF,\CG}(\CV,\CW)$ the space of maps~$\Phi\in \g_Q(\CV,\CW)$ compatible with the filtrations and define the \emph{filtered Euler number} by
\begin{align}\label{eq:general eul}
  \eul_{Q,\CF,\CG}(\CV,\CW) = \dim \g_{Q,\CF,\CG}(\CV,\CW) - \dim\CH_Q(\CV,\CW).
\end{align}
\end{definition}

\noindent
For $v\in \CH_Q(\CV)$ and $w\in \CH_Q(\CW)$, consider the map
\begin{align*}
  \delta_{v,w}\colon \g_{Q,\CF,\CG}(\CV,\CW)\to \CH_Q(\CV,\CW), \quad \Phi \mapsto \Phi v - w\Phi.
\end{align*}
Define $\Hom_{Q,\CF,\CG}(v,w) \coloneqq \ker(\delta_{v,w})$ and $\Ext_{Q,\CF,\CG}(v,w) \coloneqq \coker(\delta_{v,w})$.
Then we have an exact sequence
\begin{align*}
  0 \!\!\to\!\! \Hom_{Q,\CF,\CG}(v,w) \!\!\to\!\! \g_{Q,\CF,\CG}(\CV,\CW) \!\!\to\!\! \CH_Q(\CV,\CW) \!\!\to\!\! \Ext_{Q,\CF,\CG}(v,w) \!\!\to\!\! 0,
\end{align*}
and so \eqref{eq:general eul} is equal to
\begin{align*}
  \eul_{Q,\CF,\CG}(\CV,\CW) = \dim \Hom_{Q,\CF,\CG}(v,w)-\dim \Ext_{Q,\CF,\CG}(v,w)
\end{align*}
for any $v$ and $w$.

\begin{definition}\label{def:deltavw}
If $\eul_{Q,\CF,\CG}(\CV,\CW)=0$, then we can consider the polynomial function $P$ on $\CH_Q(\CV)\oplus \CH_Q(\CW)$ defined by
\begin{align*}
  P(v,w) \coloneqq \det(\delta_{v,w}).
\end{align*}
\end{definition}

\noindent This is a semi-invariant under the Borel subgroup~$B_Q(\CV,\CF) \times B_Q(\CW,\CG)$.

\begin{definition}
We define
$\ext_{Q,\CF,\CG}(\CV,\CW)=\min_{v,w} \dim \Ext_{Q,\CF,\CG}(v,w)$, where we minimize over all $v\in \CH_Q(\CV)$ and $w\in \CH_Q(\CW)$.
\end{definition}

The next two theorems are the analogs of Schofield's theorems~\cite{MR1162487} and the proofs are similar.

\begin{theorem}\label{theo1}
Let $(\CV,\CF)$ be a filtered dimension vector and $\CS$ a family of subspaces of~$\CV$.
Then, $\CS\subseteq_Q \CV$ if and only if
\begin{align*}
  \ext_{Q,\CF_{\CS},\CF_{\CV/\CS}}(\CS,\CV/\CS)=0.
\end{align*}
\end{theorem}
\noindent Here $\CS$ and $\CV/\CS$ are endowed with the filtrations $\CF_\CS$ and $\CF_{\CV/\CS}$ inherited from the filtrations $\CF$ on $\CV$, as defined below \cref{def:filtration}.


\begin{theorem}\label{theo2}
Let $(\CV,\CF)$ and $(\CW,\CG)$ be filtered dimension vectors.
If $\ext_{Q,\CF,\CG}(\CV,\CW)=0$, then $\eul_{Q,\CF_\CS,\CG}(\CS,\CW)\geq0$ for every $\CS\subseteq_Q\CV$.

Conversely, if $\eul_{Q,\CF_\CS,\CG}(\CS,\CW)\geq0$ for every $\CS\subseteq_Q\CV$ such that $\eul_{Q,\CF_\CS,\CF_{\CV/\CS}}(\CS,\CV/\CS)=0$, then $\ext_{Q,\CF,\CG}(\CV,\CW)=0$.
\end{theorem}

We remark that the condition~$\eul_{Q,\CF_\CS,\CF_{\CV/\CS}}(\CS,\CV/\CS)=0$ on $\CS\subseteq_Q\CV$ means that, for generic $v$, the variety of subrepresentations of~$(\CV,v)$ contained in the Schubert variety of~$\CS$ (see \cref{def:q-intersecting}) is zero-dimensional.

Now let $(\CW,\CG)$ be a filtered dimension vector and $\CV$ a family of subspaces of $\CW$.
If $\CV\subseteq_Q \CW$ then $\CU\subseteq_Q \CW$ for every $\CU\subseteq_Q \CV$;
moreover, $\eul_{Q,\CG_\CV,\CG_{\CW/\CV}}(\CV,\CW/\CV)\geq0$.
The following theorem gives, conversely, a sufficient condition for~$\CV$ to be $Q$-intersecting.

\begin{theorem}\label{theo3}
Let $\CV$ be a family of subspaces of $\CW$. Assume that $\eul_{Q,\CG_\CV,\CG_{\CW/\CV}}(\CV,\CW/\CV)\geq0$ and $\CU\subset_Q\CW$ for every proper subfamily $\CU\subset_Q \CV$.
Then, $\CV\subseteq_Q \CW$.
Moreover, it suffices to consider $\CU\subset_Q \CV$ with $\eul_{Q,\CG_\CS,\CG_{\CV/\CS}}(\CS,\CV/\CS)=0$.
\end{theorem}

The statement of this theorem is very similar to Belkale necessary and sufficient conditions for Schubert cells to intersect.
To prove \cref{theo2,theo3}, we use simplifications of Belkale's approach due to Sherman~\cite{MR3633317}, as discussed in our expository paper~\cite{BVW}.
We work directly in the context of a general quiver, which elucidates the argument.

\section{Application to the Mumford cone}\label{sec:horn applied}
Let $Q=(Q_0,Q_1)$ be a quiver and let $n=(n_x)_{x\in Q_0}$ be a dimension vector.
Consider a collection~$\CV$ of complex vector spaces $V_x$ of dimension~$n_x$.
To facilitate induction, it will be convenient not to define~$V_x=\C^{n_x}$.
Instead, we choose for each $x\in Q_0$ a subset~$J_x\subseteq\{1,2,\dots\}$ of cardinality $n_x$ and define
\begin{align*}
  V_x \coloneqq V(J_x) \coloneqq \bigoplus_{j\in J_x} \C e_j.
\end{align*}
Note that $J_x$ also naturally defines a filtration by~$F_x(j) \coloneqq \linspan \{ e_{i(1)}, \dots,$ $ e_{i(j)} \}$, where $J_x = \{ i(1) < \dots < i(n_x) \}$.
Thus, any collection $\CJ=(J_x)_{x\in Q_0}$ of finite subsets of the integers determines a collection of vectors spaces $\CV(\CJ)$ together with a collection of filtrations $\CF(\CJ)$, i.e., a filtered dimension vector.

It will be convenient to relabel all our definitions in terms of~$\CJ$.
Thus, we have the group~$\GL_Q(\CJ) \coloneqq \GL_Q(\CV(\CJ))$ acting on~$\CH_Q(\CJ) \coloneqq \CH_Q(\CV(\CJ))$, a Borel subgroup~$B_Q(\CJ)$, a unipotent subgroup~$N_Q(\CJ)$ with Lie algebra~$\n_Q(\CJ)$, and the Mumford cone~$\Cone_Q(\CJ)$.
Similarly, if~$\CJ$ and $\CL$ are two collections of subsets as above, we write $\CH_Q(\CJ,\CL) \coloneqq \CH_Q(\CV(\CJ), \CV(\CL))$.

If $\CK$ is a \emph{subfamily} of $\CJ$, i.e., a family of subsets $K_x \subseteq J_x$ for every~$x\in Q_0$, then $\CV(\CK)$ is a family of subspaces of $\CV(\CJ)$.
We will abbreviate the $Q$-intersection condition $\CV(\CK) \subseteq_Q \CV(\CJ)$ by~$\CK\subseteq_Q\CJ$ (\cref{def:q-intersecting}).
We denote
\begin{align*}
  \eul_Q(\CK,\CJ/\CK) \coloneqq \eul_Q(\CV(\CK),\CV(\CJ)/\CV(\CK)).
\end{align*}
All filtrations on subquotients of~$\CV(\CJ)$ are deduced from the filtration~$\CF(\CJ)$ of~$\CV(\CJ)$.
Here and in the following, we use the notation~$\CJ/\CK$ associated with the quotient~$\CV(\CJ)/\CV(\CK)$.

Since $K_x \subseteq J_x$, we can write $V(J_x) = V(K_x) \oplus V(J_x\setminus K_x)$ as a direct sum
and identify $V(J_x\setminus K_x) \cong V(J_x)/V(K_x)$.
Let $\n(J_x)$ denote the nilpotent radical for the $x$-th factor and define
\begin{align*}
  \n(K_x, J_x/K_x) \coloneqq \{\Phi\in \n(J_x) : \Phi(V(K_x))\subseteq V(J_x\setminus K_x), \Phi(V(J_x\setminus K_x))=0 \}.
\end{align*}
Denote by $\n(\CK,\CJ/\CK)\coloneqq\bigoplus_x \n(K_x,J_x/ K_x)$ the corresponding direct sum.
It is easy to see that
\begin{align*}
  \n(\CK,\CJ/\CK) \cong \g_{Q,\CF(\CJ)_{\CV(\CK)},\CF(\CJ)_{\CV(\CJ)/\CV(\CK)}}(\CV(\CK), \CV(\CJ)/\CV(\CK))
\end{align*}
(\cref{def:gQ}) and~\eqref{eq:general eul} amounts to
\begin{align*}
  \eul_Q(\CK,\CJ/\CK) \coloneqq \dim \n(\CK,\CJ/\CK) - \dim \CH_Q(\CK,\CJ/\CK).
\end{align*}
This is a numerical quantity that is easily computable.

\medskip

We now consider the definitions from \cref{sec:general}, with $G=\GL_Q(\CJ)$, $\g=\gl_Q(\CJ)$ and $\CH=\CH_Q(\CJ)$.
Our goal will be to describe the Mumford cone~$\Cone(G,\CH)=\Cone_Q(\CJ)$.
A basis for the Cartan subalgebra of~$\gl_Q(\CJ)$ is given by the diagonal matrices $h_{x,j}$ for $j\in J_x$ (in each factor $x\in Q_0$) such that~$h_{x,j} e_k = \delta_{j,k} e_k$ for all $k\in J_x$.
Consider~$z_x = \sum_{j\in J_x} h_{x,j}$.
Then, $z=(z_x)$ is in the center of $\gl_Q(\CV)$ and acts by zero in the infinitesimal action of $\gl_Q(\CV)$ on $\CH_Q(\CJ)$.
We can label the dominant weights for $\GL_Q(\CJ)$ by a collection $\lambda=(\lambda_x)_{x\in Q_0}$, where each $\lambda_x$ is a function $J_x\to\ZZ$ such that $\lambda_x(i) \geq \lambda_x(j)$ for all $i\leq j$.

It is easy to see that if an element $H$ is admissible (\cref{def:admissible}) then $H$ is proportional to an element of the form
\begin{align*}
  H(\CK) = (H_x), \quad H_x = \sum_{k\in K_x} h_{x,k}
\end{align*}
for some family of subsets $\CK$ of $\CJ$ (modulo possibly translations by multiples of $z=(z_x)$).
We have (\cref{def:euler mumford,def:covering ressayre}):
\begin{itemize}
\item $\eul(\CH,H)=\eul_Q(\CK,\CJ/\CK)$,
\item $H(\CK)$ is a covering element if and only if $\CK \subseteq_Q \CJ$, and therefore
\item $H(\CK)$ is a Ressayre element if and only if $\eul_Q(\CK,\CJ/\CK)=0$ and $\CK \subseteq_Q \CJ$.
\end{itemize}
As mentioned, the inequalities $\{H(\CK)\leq0\}$ for all Ressayre elements characterize the Mumford cone~\cites{ressayre2010geometric,ressayre2011geometric,VW}.
Translating this general result to our case, we obtain:

\begin{proposition}\label{prp:ressayre}
A point $\lambda=(\lambda_x)$ belongs to $\Cone_Q(\CJ)$ if and only if $\sum_{x\in Q_0,j\in J_x} \lambda_x(j) = 0$ and, for every $\CK \subseteq_Q \CJ$, we have
\begin{align*}
  \sum_{x\in Q_0} \sum_{k\in K_x} \lambda_x(k) \leq 0.
\end{align*}
Moreover, it suffices to consider $\CK$ such that $\eul_Q(\CK,\CJ/\CK)=0$.
\end{proposition}

Note that if $\CK\subseteq_Q\CJ$ then every representation~$v\in\CH_Q(\CJ)$ admits a subrepresentation with dimension vector~$\alpha_x\coloneqq\lvert K_x\rvert$.

A polynomial is semi-invariant under the full group $\GL_Q(\CJ)$ if and only if it is semi-invariant by $B_Q(\CJ)$ and its highest weight is of the form $\lambda=(\sigma_x z_x)$, i.e., it transforms as the character $\det(g)^\sigma\coloneqq\prod_x \det(g_x)^{\sigma_x}$.
It follows that the subcone~$\Sigma_Q(\CJ)$ of $\Cone_Q(\CJ)$ generated by the weights of the nonzero semi-invariant polynomials under $\GL_Q(\CJ)$ is characterized by the inequalities $\sum_x \lvert K_x\rvert \sigma_x\leq0$ for every~$\CK\subseteq_Q \CJ$.
\Cref{prp:ressayre} thus implies the description of the cone~$\Sigma_Q(\CJ)$ obtained by Derksen-Weyman~\cite{MR1758750} and Schofield-Van den Bergh~\cite{MR1908144}.
This description was also obtained using geometric invariant theory by King~\cite{MR1315461} and Ressayre~\cite{MR2875833}.
Further results were obtained by Derksen-Weyman~\cite{MR1758750} and Schofield-Van den Bergh~\cite{MR1908144} on the description of semi-invariants (under the full group $\GL_Q(\CJ))$).
In particular, they showed that the cones~$\Sigma_Q(\CJ)$ are saturated.

\medskip

We can now rephrase the results of the preceding section and characterize the $Q$-intersecting subsets $\CK\subseteq_Q \CJ$ inductively as follows.

\begin{definition}\label{def:horn set}
We define $\Horn_Q(\CJ)$ as the set of subfamilies~$\CK$ of $\CJ$ such that, if $\CK\neq\CJ$, $\eul_Q(\CK,\CJ/\CK)\geq 0$ and, for every $\CL\in \Horn_Q(\CK)$ that satisfies $\CL\neq\CK$ and $\eul_Q(\CL,\CK/\CL)=0$, we have $\eul_Q(\CL,\CJ/\CL)\geq 0$.
\end{definition}

This gives a recursive series of numerical conditions.

\begin{theorem}\label{thm:main}
Let $\CJ=(J_x)_{x\in Q_0}$ be a family of finite subsets of $\N$.  
\begin{enumerate}
\item $\CK \subseteq_Q \CJ$ if and only if $\CK\in\Horn_Q(\CJ)$.
\item A point $\lambda$ is in $\Cone_Q(\CJ)$ if and only if $\sum_{x\in Q_0,j\in J_x} \lambda_x(j) = 0$ and, for every $\CK\in\Horn_Q(\CJ)$, we have
\begin{align}\label{eq:horn ieq}
  \sum_{x\in Q_0} \sum_{k\in K_x} \lambda_x(k)\leq0.
\end{align}
It suffices to consider $\CK$ such that $\eul_Q(\CK,\CJ/\CK)=0$. 
\end{enumerate}
\end{theorem}

\begin{remark}
\Cref{thm:main} gives a complete (but in general redundant) system of inequalities for~$\Cone_Q(\CJ)$ that is straightforward to compute.
\end{remark}

Finally, let us compare \cref{def:horn set} with our general definition of a Horn element.
Thus, assume that~$H$ is a Horn element in the sense of \cref{def:horn}.
As $H$ is admissible, we may assume without loss of generality that $H=H(\CK)$ for some family~$\CK$ of subsets of $\CJ$.
We will show that $\CK\in\Horn_Q(\CJ)$.

The centralizer $G(H=0)$ of $H$ is naturally the product of two groups, $\GL_Q(\CK)$ and $\GL(\CJ\setminus\CK)$, since $H_x$ has only two eigenvalues $1$ and $0$ with respective multiplicities $K_x$ and $J_x\setminus K_x$.
Thus the Mumford cone of $\CH(H=0)$ factors accordingly as the product of the cones $\Cone_Q(\CK)$ and $\Cone_Q(\CJ\setminus\CK)$.
The Cartan subalgebra of $G(H=0)$ is a direct sum of the Cartan subalgebra of $\gl_Q(\CK)$ and of $\gl_Q(\CJ\setminus\CK)$, respectively, so this means that $\kappa_H = \kappa_{\CK} \oplus \kappa_{\CJ\setminus\CK}$, with $\kappa_{\CK} \in\Cone_Q(\CK)$ and $\kappa_{\CJ\setminus\CK}\in\Cone_Q(\CJ\setminus\CK)$.
Now, \cref{thm:main} asserts that the cone $\Cone_Q(\CK)$ is characterized by the inequalities $(H(\CL),\kappa_\CK)\leq0$ for all~$\CL\in\Horn_Q(\CK)$ such that $\eul_Q(\CL,\CK/\CL)=0$.
Under the latter hypothesis, it is easy to verify that $(H(\CL),\kappa_\CK)=-\eul_Q(\CL,\CJ/\CL)$, so it follows that, by definition, $\CK\in\Horn_Q(\CJ)$.

The preceding also implies that any Horn element is a Ressayre element.
Indeed, $\CK \subseteq_Q \CJ$ by \cref{thm:main}, so $H$ is a covering element.
But since~$H$ is a Horn element, we also have that $\eul_Q(\CK,\CJ/\CK)=0$, so $H$ is in fact a Ressayre element.
This shows that the representation~$\CH_Q(\CJ)$ has the Horn property (\cref{def:horn prop}), as we claimed in the introduction.

\end{document}